\documentclass[11pt]{elsarticle}

\usepackage{amssymb,amsmath}
\usepackage{verbatim,a4wide}
\usepackage[cp1250]{inputenc}
\usepackage[OT4]{fontenc}
\usepackage[english]{babel}

\usepackage{bm}
\usepackage{graphicx}
\usepackage{arydshln}
\usepackage{amsfonts}
\usepackage{epstopdf}
\usepackage{caption}
\usepackage{subfig}
\usepackage{appendix}
\usepackage{grffile}
\usepackage{bbm}
\usepackage{indentfirst}
\usepackage{booktabs}
\usepackage{multirow}
\usepackage{rotating}
\usepackage{hyperref}
\usepackage{color}

\numberwithin{equation}{section}
\newtheorem{thm}{Theorem}[section]
\newtheorem{lem}[thm]{Lemma}

\newtheorem{alg}[thm]{Algorithm}
\newtheorem{exmp}[thm]{Example}
\newtheorem{prob}[thm]{Problem}
\newdefinition{rem}[thm]{\it Remark }
\newproof{pf}{Proof}
\newcommand{\ra}[1]{\renewcommand{\arraystretch}{#1}}
\newcommand{\N}{\mathbb N}
\newcommand{\R}{\mathbb R}

\newcommand{\ip}[2]{\mbox{$\left\langle{#1},{#2}\right\rangle$}}
\newcommand{\der}{\mbox{${\rm\,d}$}}
\newcommand{\dt}{\mbox{${\rm\,d}t$}}
\newcommand{\du}{\mbox{${\rm\,d}u$}}

\newcommand{\rbinom}[2]{\mbox{$\displaystyle\binom{#1}{#2}^{\!\!-1}\!\!$}}

\journal{Applied Mathematics and Computation}

\begin{document}
\thispagestyle{empty}

\begin{frontmatter}

%%%%%%%%%%
\title{Efficient merging of multiple segments of B\'ezier curves}
%%%%%%%%%%

\author{Pawe{\l} Wo\'{z}ny\corref{cor}}
\ead{Pawel.Wozny@ii.uni.wroc.pl}
\author{Przemys{\l}aw Gospodarczyk}
\ead{pgo@ii.uni.wroc.pl}
\author{Stanis{\l}aw Lewanowicz}
\ead{Stanislaw.Lewanowicz@ii.uni.wroc.pl}
\cortext[cor]{Corresponding author. Fax {+}48 71 3757801}
\address{Institute of Computer Science, University of Wroc{\l}aw,
         ul.~Joliot-Curie 15, 50-383 Wroc{\l}aw, Poland}

\begin{abstract}
This paper deals with the merging problem of segments of a composite B\'ezier curve,
with the endpoints continuity constraints.
We present a novel method which is based on the idea of using constrained dual Bernstein polynomial basis
(P. Wo\'zny,  S. Lewanowicz, Comput. Aided Geom. Design 26 (2009), 566--579)
to compute the control points of the merged curve. Thanks
to  using fast schemes of evaluation of  certain connections
involving  Bernstein and dual
Bernstein polynomials,  the complexity of our algorithm
is significantly less than complexity of other merging methods.
\end{abstract}


\end{frontmatter}
%%\end{comment}
%%%%%%%%%%%%%%%%%%%%%%%%%%%%%%%%%%%%%%%%%%%%%%%%%%%%%%%%%%%%%%%%%%%%%%%%%%%
%%%%%%%%%%%%%%%%%%%%%%%%%%%%%%%%%%%%%%%%%%%%%%%%%%%%%%%%%%%%%%%%%%%%%%%%%%%

\section{Introduction}
                                                        \label{SS:Intr}
%%%%%%%%%%%%%%%%%%%%%%%%%%%%%%%%%%%%%%%%%%%%%%%%%%%%%%%%%%%%%%%%%%%%%%%%%%%
%%%%%%%%%%%%%%%%%%%%%%%%%%%%%%%%%%%%%%%%%%%%%%%%%%%%%%%%%%%%%%%%%%%%%%%%%%%

This paper deals with the merging problem of segments of a composite B\'ezier curve,
in other words:  multiple adjacent B\'ezier curves, with the endpoints continuity constraints.
More specifically, we consider the following approximation problem.
\begin{prob} \,[\textsf{Merging of multiple segments of B\'ezier curves}]
                                                \label{P:main}%PPPPPPPPPPPPPPPPPPPPPPPPPPPPP
Let $0=t_0<t_1<\ldots<t_s=1$ be a partition of the interval $[0,\,1]$.
Let be given a  composite B\'ezier curve $P(t)$  ($t\in[0,\,1]$) which in the interval $[t_{i-1},\,t_i]$ ($i=1,2,\ldots,s$)
reduces to a B\'ezier curve $P^i(t)$ of degree $n_i$, i.e.,
\begin{equation}
	\label{E:P}%EEEEEEEEEEEEEEEE
 P(t)=P^i(t):=\sum_{j=0}^{n_i}p^i_j\,B^{n_i}_j\left(\frac{t-t_{i-1}}{\Delta t_{i-1}}\right)
 \qquad (t_{i-1}\le t\le t_i),
\end{equation}
where $\Delta t_{i-1}:=t_i-t_{i-1}$, and
\[
   B^n_{j}(t):=\binom nj t^j(1-t)^{n-j} \qquad (0\le j\le n)
\]
are {Bernstein basis polynomials} of degree $n$.
Find a degree $m$ ($\ge\max_i{n_i}$)  B\'ezier curve
\begin{equation}
	\label{E:R}%EEEEEEEEEEEEEEEE
R(t):=\sum_{j=0}^m r_j\,B^m_j(t)\qquad (0\le t\le 1)
\end{equation}
such that the  error
\[	
\int_{0}^{1}\|P(t)-R(t)\|^2\dt
\]
is minimized in the space $\Pi^d_m$ of parametric polynomials in $\R^d$ of degree at most $m$
(for simplicity, we write $\Pi_m:=\Pi^1_m$)
under the additional conditions that
\begin{equation}
	\label{E:constr}%EEEEEEEEEEEEEEEE	
	\begin{array}{l}
		R^{(i)}(0)=P^{(i)}(0) \qquad (i=0,1,\ldots,k-1),\\[0.5ex]
		R^{(j)}(1)=P^{(j)}(1) \qquad (j=0,1,\ldots,l-1),
	\end{array}
\end{equation}
where $k\leq n_1+1$, $l\leq n_s+1$, and  $k+l\leq m$. Here $\|\cdot\|$ is the Euclidean vector norm.
\end{prob}
There have been many papers relevant to this problem. As for merging of two B\'ezier curves,
besides the pioneering work by Hoschek \cite{Hos87},
we should mention papers \cite{HTJS01,Lu14,Lu13,THH03,ZW09}.
Solving problem of merging more than two segments may be reduced to repeated merging of two curves.
This, however, may generate loss in accuracy of results and increase of computational cost.
The only existing algorithms to solve the problem of merging multiple B\'ezier adjacent curves
are those of \cite{CW08} and \cite{Lu15}. In the first one, only $C^0$ continuity at the endpoints
can be imposed, which results in its limited applicability in CAGD. The second algorithm
is much more general, accepting $C^{r,s}$ ($r,s\ge0$) continuity conditions. Notice that the $G^1$ multiwise
merging also was studied in \cite{Lu15}.

We present a novel method which is based on the idea of using constrained dual Bernstein polynomial basis
\cite{WL09} to compute the control points $r_i$.  Thanks
to  using  fast schemes of  evaluation  of some connections involving Bernstein and dual
Bernstein polynomials, our algorithm is rather efficient. Its  complexity is $O(sm^2)$,
which is significantly less than complexity of the methods  in \cite{CW08} and \cite{Lu15}.

The outline of this paper is as follows. Section~\ref{S:prelim} has preliminary character.
 Section~\ref{S:main} brings a complete solution to Problem~\ref{P:main}.
Section~\ref{S:Alg} deals with algorithmic implementation of the proposed method.
In Section~\ref{S:Exmp}, we give some examples showing efficiency of our method.
Conclusions are given in Section~\ref{S:Concl}.

%%%%%%%%%%%%%%%%%%%%%%%%%%%%%%%%%%%%%%%%%%%%%%%%%%%%%%%%%%%%%%%%%%%%%%%%%%%%%%%%
%%%%%%%%%%%%%%%%%%%%%%%%%%%%%%%%%%%%%%%%%%%%%%%%%%%%%%%%%%%%%%%%%%%%%%%%%%%%%%%%
\section{Preliminaries}
                                                                                                                                 \label{S:prelim}
%%%%%%%%%%%%%%%%%%%%%%%%%%%%%%%%%%%%%%%%%%%%%%%%%%%%%%%%%%%%%%%%%%%%%%%%%%%%%%%%
%%%%%%%%%%%%%%%%%%%%%%%%%%%%%%%%%%%%%%%%%%%%%%%%%%%%%%%%%%%%%%%%%%%%%%%%%%%%%%%%

Let $\Pi_m^{(k,l)}$, where $k$ and $l$ are nonnegative integers such that $k+l\le m$,
be the space of all polynomials of  degree at most $m$,
whose derivatives of order less than $k$ at $t=0$, as well as derivatives of order
less than $l$ at $t=1$, vanish:
\[
\Pi_m^{(k,l)}:=
       \left\{P\in\Pi_m\::\:
       P^{(i)}(0)=0\quad (0\le i\le k-1)
       \; \mbox{and}\;
       P^{(j)}(1)=0\quad (0\le j\le l-1)\right\}.
\]
 Obviously,  $\mbox{dim}\;\Pi_m^{(k,l)}=m-k-l+1$, and the Bernstein polynomials
\(
\left\{B^m_k,B^m_{k+1},\ldots,B^m_{m-l}\right\}
\)
form a basis of this space.
There is a unique \textit{dual constrained Bernstein basis of degree} $m$ (see, e.g., \cite{Jue98}),
\[
D^{(m,k,l)}_k,D^{(m,k,l)}_{k+1},\ldots,D^{(m,k,l)}_{m-l},
\]
satisfying
\[
	\left\langle D^{(m,k,l)}_i,\,B^m_j \right\rangle=\delta_{ij}\qquad (i,j=k,k+1,\ldots,m-l),
\]
where $\delta_{ij}$ is 1 if $i=j$ and 0 otherwise, and the inner product $\langle\cdot,\cdot\rangle$ is given by
\[
	\langle f,\,g \rangle:=\int_{0}^{1} f(t)g(t)\dt.
\]
For $k=l=0$ (the unconstrained case), we have
dual Bernstein basis $D^m_i:=D^{(m,0,0)}_i$
($i=0,1,\ldots,m$) of the space $\Pi_m^{(0,0)}=\Pi_m$.

\begin{lem}\label{L:BernToDual} %%LLLLLLLLLLLLLLLLLLLLLLLLL
Let $n$ and $m$ be positive integers such that $n \leq m$. The following formula holds:
\[
B_i^n(t) = \sum_{j=0}^m a_{ij}^{(n,m)}D^m_j(t) \qquad (0 \leq i \leq n;\ n \leq m),
\]
where
\begin{equation}\label{E:CoeffBernToDual}
a_{ij}^{(n,m)} := \frac{1}{m+n+1}\binom{n}{i}\binom{m}{j}\rbinom{n+m}{i+j}.
\end{equation}
\end{lem}
\begin{pf}
Obviously, we have
\[
	a_{ij}^{(n,m)}=\langle B_i^n,B_j^m\rangle=\int_{0}^{1}B^{n}_i(t) B^{m}_j(t)\dt,
\]
and the result follows by  the well known properties of Bernstein polynomials (see, e.g., \cite[\S6.10]{Far02}):
\begin{align*}
&B^{n}_i(t) B^{m}_j(t)=\binom{n}{i}\binom{m}{j}\rbinom{n+m}{i+j}B^{n+m}_{i+j}(t),\\
 &\int_{0}^{1}B^{n+m}_{i+j}(t)\dt=\frac1{n+m+1}.	
\end{align*}
\qed
\end{pf}
\begin{lem}\label{L:bestpol}%%LLLLLLLLLLLLLLLLLLLLLLLLLLLLL
Let $m,k,l\in\N$ be such that $0\le k+l\le m$ and let  $f$ be a function defined
on $[0,\,1]$.
The polynomial $S\in\Pi_m^{(k,l)}$,  which gives minimum
value of the norm
\[
	 \|f-S\|_{L_2}:=\ip{f-S}{f-S}^{\frac12},
\]
is given by
\begin{equation}
	\label{E:bestpol}%EEEEEEEEEEEEEEEE
	S=\sum_{i=k}^{m-l}\ip f{B^m_i}D^{(m,k,l)}_i.
\end{equation}
\end{lem}
\begin{pf}
	Obviously, $S$ has the following representation in the dual Bernstein basis of the space
	$\Pi_m^{(k,l)}$:
	\[
		S=\sum_{i=k}^{m-l}\ip {S}{B^m_i} D^{(m,k,l)}_i.
	\]
On the other hand, a classical characterization of the best approximation polynomial $S$ is that $\ip {f-S}Q=0$ holds for
any polynomial  $Q\in\Pi_m^{(k,l)}$. In particular, for $Q=B^m_i$, we obtain
\[
	\ip {f}{B^m_i}=\ip {S}{B^m_i} \qquad (k\le i\le m-l).
\]
Hence, the formula \eqref{E:bestpol} follows.
\qed\end{pf}

Further properties of  the  polynomials $D^{(m,k,l)}_i$ are studied in  \cite{LW11a,WL09} and in
the recent paper \cite{LW11b}, where the following result is given.
\begin{lem}[\cite{LW11b}] \label{L:constrDinB}%LLLLLLLLLLLLLLLLLLLLLLLLL
The constrained dual basis polynomials have the  B\'ezier-Bernstein representation
\begin{equation}\label{E:constrDinB}%EEEEEEEEEEEEEEEE
	 D^{(m,k,l)}_i=\sum_{j=k}^{m-l}c_{ij}(m,k,l)\,
	                   B^{m}_j,
\end{equation}
where the coefficients $c_{ij}\equiv c_{ij}(m,k,l)$
satisfy the recurrence relation
\begin{align}\label{E:c-rec}%%EEEEEEEEEEEEEEEEEEEEEEEEEEEEE
c_{i+1,j}=&\frac1{A(i)}\,\left\{2(i-j)(i+j-m)\,c_{ij}
+B(j)\,c_{i,j-1} +A(j)\,c_{i,j+1}-B(i)\,c_{i-1,j}\right\}\nonumber \\[1.25ex]
&\hphantom{\frac1{A(i)}\,\{+B(j)\,c_{i,j-1} +A(j)\,c_{i,j+1}}
(k\le i\le m-l-1,\quad  k\le j\le m-l)
\end{align}	
with
\[
\begin{array}{l}
A(u):=(u-m)(u-k+1)(u+k+1)/(u+1),\\[2ex]	
B(u):=u(u-m-l-1)(u-m+l-1)/(u-m-1).
\end{array}
\]
We adopt the convention that
$c_{ij}:=0$ if $i<k$, or $i>m-l$, or $j<k$, or $j>m-l$.
The starting values are
\begin{equation}
	\label{E:c-start}%EEEEEEEEEEEEEEEE
c_{kj}
 =(-1)^{j-k}(2k+1)\rbinom{m}{k}\,\binom{m+k-l+1}{2k+1}
 \rbinom{m}{j}\,\binom{m-k-l}{j-k}\binom{m+k+l+1}{k+j+1}, 	
\end{equation}
where $j=k,k+1,\ldots, m-l$.
\end{lem}

In the next section, we will need the following restriction of the representation
of the polynomial $B^m_j$ to a subinterval  of the interval $[0,\,1]$.
 \begin{lem}							\label{C:subdivB}%%LLLLLLLLLLLLLLLLLLL
 	 Let $0=t_0<t_1<\ldots<t_s=1$ be a partition of the interval $[0,\,1]$.
In the subinterval $[t_{i-1},\,t_i]$ ($i=1,2,\ldots,s$), the Bernstein polynomial $B^m_j$
can be expressed in the form
\begin{equation}
	\label{E:subdivB}%EEEEEEEEEEEEEEEE
	B^m_j(t)=\sum_{h=0}^{m}d^{(i)}_{jh}B^m_h\left(\frac{t-t_{i-1}}{\Delta t_{i-1}}\right),
\end{equation}
where
\begin{equation}
	\label{E:dijh}%EEEEEEEEEEEEEEEE
	d^{(i)}_{jh}:=\sum_{v=0}^{h}B^{m-h}_{j-v}(t_{i-1})B^h_v(t_i).
\end{equation}
\end{lem}
\begin{pf}
The result is obtained in two steps. First, subdivide the polynomial
\[
	B^m_j(t)=\sum_{h=0}^{m}\delta_{jh}B^m_h(t)
\]
at the point $t_i$ to obtain two forms for
the subintervals $[0,\,t_i]$ and $[t_i,\,1]$. Next, subdivide the form corresponding to $[0,\,t_i]$ at $t_{i-1}/t_i$.
We obtain the formula \eqref{E:subdivB} with the coefficients $d^{(i)}_{jh}$ given by
\[
	d^{(i)}_{jh}:=\sum_{w=0}^{m-h}B^{m-h}_{w}(t_{i-1}/t_i)B^{w+h}_j(t_i)
\]
(we ignore the fact that the initial terms of the sum vanish as  $B^{w+h}_j(t_i)=0$ for $0\le w<j-h$).
Using the identity
\[
	 B^{n+q}_j(x)=\sum_{w=0}^{q}B^q_w(x)B^n_{j-w}(x),
\]
which can be easily proved using some basic properties of the Bernstein polynomials (see, e.g., \cite[\S6.10]{Far02}),
and
\[
	B^n_j(cx)  = \sum_{v=0}^{n}B^n_v(x)B^v_j(c)
\]
(\textit{ibid.}), it can be seen that
\begin{align*}
	  d^{(i)}_{jh}&=  \sum_{w=0}^{m-h}B^{m-h}_{w}(t_{i-1}/t_i)\sum_{v=0}^{h}B^{h}_v(t_i)B^w_{j-v}(t_i) \\
	  &=\sum_{v=0}^{h}B^{h}_v(t_i)\sum_{w=0}^{m-h}B^{m-h}_{w}(t_{i-1}/t_i)B^w_{j-v}(t_i)\\
	  &=\sum_{v=0}^{h}B^{h}_v(t_i)B^{m-h}_{j-v}(t_{i-1}).
\end{align*}
\qed\end{pf}
Equation \eqref{E:subdivB} is obviously equivalent to
\begin{equation}
	\label{E:subdivB-var}%EEEEEEEEEEEEEEEE	
	B^m_j(u\Delta t_{i-1}+t_{i-1})=\sum_{h=0}^{m}d^{(i)}_{jh}B^m_h\left(u\right)\qquad (0\le u\le1).
\end{equation}
Now, by the bi-orthogonality property of the bases $\{B^m_h\}$ and $\{D^m_g\}$, we have
\begin{equation}
	\label{E:dijh-var}%EEEEEEEEEEEEEEEE
	d^{(i)}_{jh}=\int_{0}^{1}B^m_j(u\Delta t_{i-1}+t_{i-1})D^m_h(u)\du.
\end{equation}

\begin{lem}	                                    \label{L:RecCoeffsD}%LLLLLLLLLLLLLLLLLLLLLLLLLLLLLLLLL
For $i=1,2,\ldots s$, the coefficients  $d_{jh}^{(i)}$ satisfy the following recurrence equation:
\begin{multline*}
 \Delta t_{i-1}\left[(m-j+1)d_{j-1,h}^{(i)}+(2j-m)d_{jh}^{(i)}-(j+1)d_{j+1,h}^{(i)}\right]\\
\qquad\quad = (m-h)d_{j,h+1}^{(i)} + (2h-m)d_{jh}^{(i)} - hd_{j,h-1}^{(i)}\\
 (1 \leq j \leq m-1;\ 0 \leq h \leq m).
\end{multline*}
\end{lem}
\begin{pf}
Differentiate both sides of Equation \eqref{E:subdivB-var} with respect to $u$, and make use of the identity
\[
\frac{\mbox{d}}{\mbox{d}u}B^m_j(u) = (m-j+1)B^m_{j-1}(u)+(2j-m)B^m_j(u)-(j+1)B^m_{j+1}(u).
\]
Equating the  B\'ezier coefficients gives the result.
\qed
\end{pf}
%%%%%%%%%%%%%%%%%%%%%%%%%%%%%%%%%%%%%%%%%%%%%%%%%%%%%%%%%%%%%%%%%%%%%%%%%%%%%%%%%
%%%%%%%%%%%%%%%%%%%%%%%%%%%%%%%%%%%%%%%%%%%%%%%%%%%%%%%%%%%%%%%%%%%%%%%%%%%%%%%%%

%%%%%%%%%%%%%%%%%%%%%%%%%%%%%%%%%%%%%%%%%%%%%%%%%%%%%%%%%%%%%%%%%%%%%%%%%%%%%%%%%
%%%%%%%%%%%%%%%%%%%%%%%%%%%%%%%%%%%%%%%%%%%%%%%%%%%%%%%%%%%%%%%%%%%%%%%%%%%%%%%%%
\section{Merging of the composite B\'ezier curve  segments}
						\label{S:main}%SSSSSSSSSSSSSSSSSSSSSSSSSS
%%%%%%%%%%%%%%%%%%%%%%%%%%%%%%%%%%%%%%%%%%%%%%%%%%%%%%%%%%%%%%%%%%%%%%%%%%%%%%%%%
%%%%%%%%%%%%%%%%%%%%%%%%%%%%%%%%%%%%%%%%%%%%%%%%%%%%%%%%%%%%%%%%%%%%%%%%%%%%%%%%%
Clearly, the B\'ezier curve being the solution of Problem~\ref{P:main} can be obtained
in a componentwise way. Hence, it is sufficient to give the details of our method of solving
this problem in case where $d=1$.

%%%%%%%%%%%%%%%%%%%%%%%%%%%%%%%%%%%%%%%%%%%%%%%%%%%%%
\begin{thm}
	\label{T:main}%TTTTTTTTTTTTTTTT
Let $0=t_0<t_1<\ldots<t_s=1$ be a partition of the interval $[0,\,1]$.
Let be given the piecewise polynomial function $P(t)$  ($t\in[0,\,1]$),
which in the interval $[t_{i-1},\,t_i]$ ($i=1,2,\ldots,s$)
reduces to a polynomial  $P^i(t)$ of degree $n_i$,	
with the B\'ezier coefficients $p^i_j$ ($i=1,2,\ldots,s;\; j=0,1,\ldots,n_i$)
(cf. \eqref{E:P}). The coefficients $r_0,r_1, \ldots,r_m$ of the polynomial
\eqref{E:R} minimising the error
\[
	\|R-P\|^2_{L_2}:=\langle R-P,R-P\rangle
\]
with constraints \eqref{E:constr}
are given by
\begin{align}
         			\label{E:ri-begin}%EEEEEEEEEEEEEEEE
	  & \displaystyle
	  r_j=\binom{n_1}{j}\binom{m}{j}^{\!-1}\,\Delta^jp^1_0-
          \sum_{h=0}^{j-1}(-1)^{j+h}\binom{j}{ h}r_{h}\qquad (j=0,1,\ldots,k-1),\\[1ex]
	  			\label{E:ri-end}%EEEEEEEEEEEEEEEEE
	& \displaystyle r_{m-j}
	=(-1)^j\binom{n_s}{j}\binom{m}{j}^{\!-1}\,\Delta^jp^s_{n_s-j}-
	\sum_{h=1}^{j}(-1)^h\binom{j}{ h}r_{m-j+h}\qquad (j=0,1,\ldots,l-1),	 \\[1ex]	
	 \label{E:ri-mid}%EEEEEEEEEEEEEEEE
	  &r_{j}= \sum_{h=k}^{m-l}\hat{r}_hc_{hj}(m,k,l)\qquad (j=k,k+1,\ldots,m-l), 	  	
\end{align}
where
\begin{align}
	\label{E:hatr}%EEEEEEEEEEEEEEEE
	  \hat{r}_h&:=\sum_{i=1}^{s}\Delta t_{i-1}\sum_{v=0}^{m}\hat{p}^i_vd^{(i)}_{hv}
	  -\frac1{2m+1}\binom{m}{h}\left(\sum_{v=0}^{k-1}+\sum_{v=m-l+1}^{m}\right)\rbinom{2m}{h+v}\binom{m}{v}r_v,  \\
	\label{E:hatp}%EEEEEEEEEEEEEEEE
	  \hat{p}^i_v&:= \frac{1}{m+n_i+1}\binom{m}{v}\sum_{q=0}^{n_i}\rbinom{m+n_i}{q+v}\binom{n_i}{q}p^i_q,
\end{align}
with $c_{hj}(m,k,l)$ and $d^{(i)}_{jh}$ being introduced in \eqref{E:constrDinB} and \eqref{E:dijh}, respectively.
Here we use the standard notation $\Delta^0c_h:=c_h$, $\Delta^j
c_h:=\Delta^{j-1}c_{h+1}-\Delta^{j-1}c_{h}$ ($j=1,2,\ldots$).
\end{thm}

\begin{pf} Recall that for arbitrary polynomial of degree $N$,
\[
	U_N(t)=\sum_{h=0}^{N}u_h\,B^N_h(t),
\]
the well-known formulas hold (see, e.g., \cite[\S5.3]{Far02})
\begin{eqnarray*}
U^{(j)}_N(0)&=&\frac{N!}{(N-j)!}\Delta^ju_0=
\frac{N!}{(N-j)!}\sum_{h=0}^{j}(-1)^{j+h}\binom{j}{h}u_{h},\\
U^{(j)}_N(1)&=&\frac{N!}{(N-j)!}\Delta^j u_{N-j}=
\frac{N!}{(N-j)!}\sum_{h=0}^{j}(-1)^{j+h}\binom{j}{h}u_{N-j+h}.
\end{eqnarray*}
Using the above equations
in \eqref{E:constr}, we obtain the forms
\eqref{E:ri-begin} and
\eqref{E:ri-end} for the coefficients $r_0,\,r_1,\ldots,r_{k-1}$ and
$r_{m-l+1},\ldots,r_{m-1},r_m$, respectively.

The remaining coefficients $r_k,\,r_{k+1},\,\ldots,\,r_{m-l}$ are to be determined so that
\[
\| P-R\|^2_{L_2}=\|W-S\|^2_{L_2}	
\]
has the least value, where
\begin{align*}
W&:=P-\left(\sum_{h=0}^{k-1}+\sum_{h=m-l+1}^{m}\right)r_hB^m_h,\\
S&:=\sum_{j=k}^{m-l}r_jB^m_j.
\end{align*}
To be strict, we first obtain the coefficients $\hat{r}_j$  of the searched polynomial
in the constrained dual  Bernstein basis $\{D^{(m,k,l)}_h\}$,
\[
	S=\sum_{j=k}^{m-l}\hat{r}_jD^{(m,k,l)}_j;
\]
then the B\'ezier coefficients $r_j$ of $S$ will be easily computed using Equation \eqref{E:ri-mid} (cf. Lemma~\ref{L:constrDinB}).

Now, using Lemma~\ref{L:BernToDual}, we represent each segment $P^i$ of the original piecewise polynomial $P$
in the dual  Bernstein basis of degree $m$,
\[
P^i(t) %%= \sum_{h=0}^{n_i}p^i_hB^{n_i}_h\left(\frac{t-t_{i-1}}{\Delta t_{i-1}}\right)
 = \sum_{v=0}^m \hat{p}^i_vD^{m}_v\left(\frac{t-t_{i-1}}{\Delta t_{i-1}}\right)
\]
with $\hat{p}^i_v$ being defined in \eqref{E:hatp}.

Using Lemma~\ref{L:bestpol} and Equation \eqref{E:dijh-var}, we obtain
\begin{align*}
\hat{r}_j = & \left<W,B_j^m\right>  = \int_0^1W(t)B_j^m(t)\dt\\
	=&\sum_{i=1}^{s}\sum_{h=0}^m\hat{p}^i_h
              \int_{t_{i-1}}^{t_{i}}D_h^m\left(\frac{t-t_{i-1}}{\Delta t_{i-1}}\right)B_j^m(t)\dt\\
            &\hphantom{sum_{i=1}^{s}} - \left(\sum_{h=0}^{k-1}+\sum_{h=m-l+1}^{m}\right)r_h\int_0^1B_h^m(t)B_j^m(t)\dt\\
             =&
              \sum_{i=1}^{s}\sum_{h=0}^m\hat{p}^i_h\Delta t_{i-1}
              \int_0^1D_h^m(u)B_j^m(\Delta t_{i-1}u+t_{i-1})\du \\
             &\hphantom{sum_{i=1}^{s}}- \left(\sum_{h=0}^{k-1}+\sum_{h=m-l+1}^{m}\right)
              r_h \frac{1}{2m+1}\binom{m}{h}\binom{m}{j}\rbinom{2m}{h+j}\\
              =&\sum_{i=1}^{s}\Delta t_{i-1}\sum_{h=0}^m\hat{p}^i_h d_{jh}^{(i)}
              %\\ &\hphantom{sum_{i=1}^{s}}
              - \frac1{2m+1}\binom{m}{j}\left(\sum_{h=0}^{k-1}+\sum_{h=m-l+1}^{m}\right)
              r_h \binom{m}{h}\rbinom{2m}{h+j}\\[1ex]
              &\hphantom{sum_{i=1}^{s}\sum_{h=0}^m\hat{p}^i_h\int_{t_{i-1}}^{t_{i}}D_h^m\left(\frac{t-t_{i-1}}{\Delta t_{i-1}}\right)B_j^m(t)\dt}
              (j=k,k+1,\ldots,m-l).
\end{align*}
This completes the proof.
\qed\end{pf}

Now, let the composite curve $P$ and the merged curve $R$ be curves in $\mathbb{R}^d$ ($d\ge1$).
Let $p^i_{j}=(p^i_{j1},p^i_{j2},\ldots,p^i_{jd})$ ($i=1,2,\ldots,s;\ j=0,1,\ldots,n_i$), and
$r_j=(r_{j1},r_{j2},\ldots,r_{jd})$ ($j=0,1,\ldots,m$) be the control points of  $P$
and  $R$, respectively. For $i=1,2,\ldots,s$ and $h=1,2,\ldots,d$, let us define vectors
\begin{align*}
\pi^i_h :=& \left[p^i_{0h},p^i_{1h},\ldots,p^i_{n_i,h}\right] \in \mathbb{R}^{n_i+1},\\
\varrho^i_h := &\left[\varrho^i_{0h},\varrho^i_{1h},\ldots,\varrho^i_{mh}\right] \in \mathbb{R}^{m+1},
\end{align*}
where
\begin{equation}\label{E:tilr}%%EEEEEEEEEEEEEEEEEEEEEEEEEEEE
\varrho^{i}_{zh} := \sum_{j=0}^{m}r_{jh}d^{(i)}_{jz} \qquad (z=0,1,\ldots,m).
\end{equation}
It can be shown that the $L_2$-distance between the curves $P$ and $R$ is given by the formula:
\begin{align}
E_2 :=& \| P-R\|_{L_2}\nonumber \\
      =&  \left(\sum_{i=1}^s \Delta t_{i-1} \sum_{h=1}^d\left[
                         I_{n_i,n_i}(\mathbf{\pi}^i_h,\mathbf{\pi}^i_h) - 2I_{n_i,m}(\mathbf{\pi}^i_h,\mathbf{\varrho}^i_h) + I_{mm}(\mathbf{\varrho}^i_h,\mathbf{\varrho}^i_h)
                         \right]\right)^{\frac{1}{2}},\label{E:dist}%EEEEEEEEEEEEEEEEEEEEEE
\end{align}
where
\[
I_{NM}(u,v) := \sum_{j=0}^{N}u_j\sum_{z=0}^{M} a_{jz}^{(N,M)}v_z,
\]
with $u:= \left[u_0,u_1,\ldots,u_N\right]$ and $v := \left[v_0,v_1,\ldots,v_M\right]$,
the notation used being that of \eqref{E:CoeffBernToDual}.

%%%%%%%%%%%%%%%%%%%%%%%%%%%%%%%%%%%%%%%%%%%%%%%%%%%%%%%%%%%%%%%%%%%%%%%%%%%%%
%%%%%%%%%%%%%%%%%%%%%%%%%%%%%%%%%%%%%%%%%%%%%%%%%%%%%%%%%%%%%%%%%%%%%%%%%%%%%
\section{Algorithms}
                     \label{S:Alg}%SSSSSSSSSSSSSSSSSSSSSSSSSSSSSSSSSSSSSSSS
%%%%%%%%%%%%%%%%%%%%%%%%%%%%%%%%%%%%%%%%%%%%%%%%%%%%%%%%%%%%%%%%%%%%%%%%%%%%%
%%%%%%%%%%%%%%%%%%%%%%%%%%%%%%%%%%%%%%%%%%%%%%%%%%%%%%%%%%%%%%%%%%%%%%%%%%%%%

\subsection{Auxiliary computations}
                     \label{SS:AuxAlg}%SSSSSSSSSSSSSSSSSSSSSSSSSSSSSSSSSSSSSSSS
%%%%%%%%%%%%%%%%%%%%%%%%%%%%%%%%%%%%%%%%%%%%%%%%%%%%%%%%%%%%%%%%%%%%%%%%%%%%%
%%%%%%%%%%%%%%%%%%%%%%%%%%%%%%%%%%%%%%%%%%%%%%%%%%%%%%%%%%%%%%%%%%%%%%%%%%%

In this section, we discuss details of algorithmic implementation of the results given in Theorem~\ref{T:main}.
First, we have to precompute efficiently the  coefficients  $ c_{ij}(m,k,l)\;$
introduced in Lemma~\ref{L:constrDinB} (see~Table~\ref{Tab:C}).
\begin{table}[h]
\captionsetup{margin=0pt, font={scriptsize}}
\[
\begin{array}{cccccc}
    &0&0&\ldots&0& \\
    0& c_{kk}&c_{k,k+1}&\ldots&c_{k,m-l}&0     \\
    0& c_{k+1,k}&c_{k+1,k+1}&\ldots&c_{k+1,m-l}&0  \\
    \multicolumn{6}{c}{\dotfill}\\    0& c_{m-l,k}&c_{m-l,k+1}&\ldots&c_{m-l,m-l}&0  \\
     &0&0&\ldots&0&
\end{array}
\]	
\caption{The $c$-table \label{Tab:C}}
\end{table}

Now,  the table can be completed easily by using formulas \eqref{E:c-rec}, \eqref{E:c-start}  (cf. \cite[Algorithm 3.3]{LW11b}), with the complexity $O(m^2)$.

Another task is to evaluate all the coefficients $d_{jh}^{(i)}$  ($i=1,2,\ldots,s;\ j=0,1,\ldots,m;\ h=0,1,\ldots,m$)
(cf. \eqref{E:dijh}). Thanks to Lemma~\ref{L:RecCoeffsD},  we can do it using the following algorithm.
\
\begin{alg}\![\textsf{Evaluation of the coefficients }  $d_{jh}^{(i)}$]
        \label{A:dtab}\\%AAAAAAAAAAAAAAAAAAAAAAAAAAAAAAAA
\texttt{Input}: $m$, $s$,  $0=t_0<t_1<\ldots<t_s=1$\\
\texttt{Output}: table of the coefficients $d_{jh}^{(i)}\ (i=1,2,\ldots,s;\ j=0,1,\ldots,m;\ h=0,1,\ldots,m)$
\begin{description}
\itemsep2pt
\item[\texttt{Step 1}.] For $i=1,2,\ldots,s$, compute
\begin{align*}
&d_{-10}^{(i)} := 0,\quad d_{00}^{(i)} := (1-t_{i-1})^m,\\
&d_{-1h}^{(i)} := 0,\quad
d_{0h}^{(i)}:=\frac{1-t_i}{1-t_{i-1}}d_{0,h-1}^{(i)}\qquad (h=1,2,\ldots,m).
\end{align*}
\item[\texttt{Step 2}.]   For $i=1,2,\ldots,s$,  $j=0,1,\ldots,m-1$, and $h=0,1,\ldots,m$, compute
\begin{align*}
d_{j+1,h}^{(i)} := &(j+1)^{-1}\left\{\left(\Delta t_{i-1}\right)^{-1}
	\left[hd_{j,h-1}^{(i)} - (2h-m)d_{jh}^{(i)}-(m-h)d_{j,h+1}^{(i)}\right]\right. \\
		&\hphantom{(j+1)^{-1}\{} \left. +(m-j+1)d_{j-1,h}^{(i)}+(2j-m)d_{jh}^{(i)}\right\}.
\end{align*}
\end{description}
\end{alg}

Observe that complexity of  Algorithm~\ref{A:dtab} is $O(sm^2)$.
%%%%%%%%%%%%%%%%%%%%%%%%%%%%%%%%%%%%%%%%%%%%%%%%%%%%%%%%%%%%%%%%%%%%%%%%%%%%%
%%%%%%%%%%%%%%%%%%%%%%%%%%%%%%%%%%%%%%%%%%%%%%%%%%%%%%%%%%%%%%%%%%%%%%%%%%%%%

\subsection{Main algorithm}
                     \label{SS:MainAlg}%SSSSSSSSSSSSSSSSSSSSSSSSSSSSSSSSSSSSSSSS
%%%%%%%%%%%%%%%%%%%%%%%%%%%%%%%%%%%%%%%%%%%%%%%%%%%%%%%%%%%%%%%%%%%%%%%%%%%%%
%%%%%%%%%%%%%%%%%%%%%%%%%%%%%%%%%%%%%%%%%%%%%%%%%%%%%%%%%%%%%%%%%%%%%%%%%%%%%
Now, the presented method of merging of segments of a composite B\'ezier curve is summarized in the following algorithm.
\
\begin{alg}\,
[\textsf{Merging of segments of a composite B\'ezier curve}]
	\label{A:MainAlg}\\%AAAAAAAAAAAAAAAAAAAAAAAAAA	
\texttt{Input}:  $p^i_j\ (j=0,1,\ldots,n_i)$, $n_i\ (i=1,2,\ldots,s)$, \\
\hphantom{\texttt{Input}:  }$m$, $k$, $l$, $0=t_0<t_1<\ldots<t_s=1$\\
\noindent \texttt{Output}: solution $r_0,r_1,\dots,r_m$ of the Problem~\ref{P:main}, and its error $E_2$
\begin{description}
\itemsep2pt
\item[\texttt{Step 1}.] Compute $r_0,r_1,\ldots, r_{k-1}$  by \eqref{E:ri-begin}.
\item[\texttt{Step 2}.] Compute  $r_{m-l+1},r_{m-l+2},\ldots, r_m$  by \eqref{E:ri-end}.
\item[\texttt{Step 3}.] Compute $\hat{p}^i_j \ (i=1,2,\ldots,s;\ j=0,1,\ldots,m)$ by \eqref{E:hatp}.
\item[\texttt{Step 4}.] Compute $d_{jh}^{(i)}$ for $i=1,2,\ldots,s$; $j=0,1,\ldots,m$; $h=0,1,\ldots,m$,
			using Algorithm~\ref{A:dtab}.
\item[\texttt{Step 5}.] Compute $\hat{r}_j\ (j=k,k+1,\ldots,m-l)$  by \eqref{E:hatr}.
\item[\texttt{Step 6}.] Compute $c_{ij}(m,k,l)$ for $i,j=k,k+1,\ldots,m-l$, using \eqref{E:c-rec}, \eqref{E:c-start}  (cf. \cite[Algorithm 3.3]{LW11b}).
\item[\texttt{Step 7}.] Compute $r_j\ (j=k,k+1,\ldots,m-l)$ by \eqref{E:ri-mid}.
\item[\texttt{Step 8}.] Compute $\varrho^{i}_{zh}\ (i=1,2,\ldots,s;\ z=0,1,\ldots,m;\ h=1,2,\ldots,d)$ by \eqref{E:tilr}.
\item[\texttt{Step 9}.] Compute $E_2$ by \eqref{E:dist}.
\end{description}
\end{alg}

Notice that complexity of Algorithm~\ref{A:MainAlg} is $O(sm^2)$.

%%%%%%%%%%%%%%%%%%%%%%%%%%%%%%%%%%%%%%%%%%%%%%%%%%%%%%%%%%%%%%%%%%%%%%%%%%%%%
%%%%%%%%%%%%%%%%%%%%%%%%%%%%%%%%%%%%%%%%%%%%%%%%%%%%%%%%%%%%%%%%%%%%%%%%%%%%%
\section{Examples}
                     \label{S:Exmp}%SSSSSSSSSSSSSSSSSSSSSSSSSSSSSSSSSSSSSSSS
%%%%%%%%%%%%%%%%%%%%%%%%%%%%%%%%%%%%%%%%%%%%%%%%%%%%%%%%%%%%%%%%%%%%%%%%%%%%%
%%%%%%%%%%%%%%%%%%%%%%%%%%%%%%%%%%%%%%%%%%%%%%%%%%%%%%%%%%%%%%%%%%%%%%%%%%%%%
In this section, we give several examples of using Algorithm~\ref{A:MainAlg}.
In every case we give the $L_2$-error $E_2$ as well as the maximum error
\[
E_{\infty} := \max_{t \in D_N} \|P(t) - R(t)\| \approx \max_{t \in [0,1]} \|P(t) - R(t)\|,
\]
where $D_N := \left\{0, 1/N, 2/N,\ldots, 1\right\}$ with $N= 500$.
Generalizing the approach of \cite[(6.1)]{Lu14}, partition of the interval $[t_0, t_s]=[0,\,1]$
is determined according to the lengths of segments $P^i$:
\begin{equation}\label{E:part}
t_j :=L_j/L_s  \qquad (j=1,2,\ldots,s-1),
\end{equation}
where
\[
L_q:=\sum_{i=1}^{q}\int_{0}^{1}\left\Vert\frac{\der}{\dt}\sum_{h=0}^{n_i}p^i_hB^{n_i}_h(t)\right\Vert\dt.
\]
Integrals are evaluated using the $\mbox{Maple}^{TM}13$ function \texttt{int} with the option \texttt{numeric}.

Results of the experiments have been obtained on a computer with \texttt{Intel Core i5-3337U 1.8GHz} processor and \texttt{8GB} of \texttt{RAM}, using $32$-digit arithmetic. Notice that $\mbox{Maple}^{TM}13$ worksheet containing programs and tests
can be found on the webpage webpage \url{http://www.ii.uni.wroc.pl/~pgo/papers.html}.

\begin{exmp}\label{Ex:Ampersand}
We use Algorithm~\ref{A:MainAlg} to merge the composite curve ``Ampersand'', with three fifth degree B\'{e}zier segments, defined by the control points $\{(1.09,0.03),$ $(1.02,0.21),$ $(0.6,0.75),$ $(0.5,1.11),$ $(0.85,1.12),$ $(0.93,1.03)\}$,
$\{(0.93,1.03),$ $(1.01,0.96),$ $(1.02,0.76),$ $(0.8,0.65),$ $(0.62,$\\$0.38),$ $(0.61,0.23)\}$, and
$\{(0.61,0.23),$ $(0.59,0.1),$ $(0.67,0.02),$ $(0.91,-0.05),$ $(1.12,0.05),$ $(1.08,$\\$0.22)\}$, respectively.
According to \eqref{E:part}, we have  $t_0 = 0,\ t_1 \doteq 0.45,\ t_2 \doteq 0.76,\ t_3 = 1$.
Obtained results are given in Table~\ref{Tab:ex1}. Moreover, we give the comparison of running times required to compute the resulting control points.
Clearly, our method is faster than the one presented in~\cite{Lu15}.
Figures~\ref{Fig:1a}~and~\ref{Fig:1b} illustrate the results
for two representative cases.
This example shows that merging may result in data compression.
\begin{table}[h]
\captionsetup{margin=0pt, font={scriptsize}}
\ra{1}
\centering
\scalebox{1}{
\begin{tabular}{@{}ccccccccc@{}}
\toprule \multicolumn{3}{c}{Parameters} & \phantom{abc} & \multicolumn{2}{c}{Errors} & \phantom{abc} & \multicolumn{2}{c}{Running times [ms]}
\\ \cmidrule{1-3} \cmidrule{5-6} \cmidrule{8-9}
$m$ & $k$ & $l$ & \phantom{abc} & $E_{2}$ & $E_{\infty}$ & \phantom{abc} & Algorithm~\ref{A:MainAlg} &  Lu~\cite{Lu15}
\\ \midrule
$8$ & $2$ & $1$ & & $8.57E{-}3$ & $2.36E{-}2$ &  & $10$ &$85$\\
  &   $2$ & $2$ & & $1.99E{-}2$ & $5.46E{-}2$ &  & $11$ & $87$\\
  &   $3$ & $2$ & & $3.89E{-}2$ & $1.04E{-}1$ &  & $10$ & $77$\\
\midrule
$10$ & $2$ & $1$ & & $3.49E{-}3$ & $1.32E{-}2$ & & $16$ & $121$\\
  &  $2$ & $2$ & & $9.43E{-}3$ & $3.36E{-}2$ & & $15$ & $108$\\
  &  $3$ & $2$ & & $1.98E{-}2$ & $6.08E{-}2$ & & $15$ & $104$\\
  \midrule
$12$ & $2$ & $1$ & & $2.70E{-}3$ & $9.84E{-}3$ & & $22$ & $167$\\
  & $2$ & $2$ & & $5.71E{-}3$ & $2.29E{-}2$ & & $22$ & $160$\\
  &   $3$ & $2$ & & $1.06E{-}2$ & $3.81E{-}2$ & & $19$ & $158$\\
\end{tabular}}
\caption{Least-squares and maximum errors for merging of three  segments of the composite B\'{e}zier curve with constraints.}
\label{Tab:ex1}
\end{table}
\newpage
\begin{figure}[h]
\captionsetup{margin=0pt, font={scriptsize}}
\begin{center}
\setlength{\tabcolsep}{0mm}
\begin{tabular}{c}
\subfloat[]{\label{Fig:1a}\includegraphics[width=0.38\textwidth]{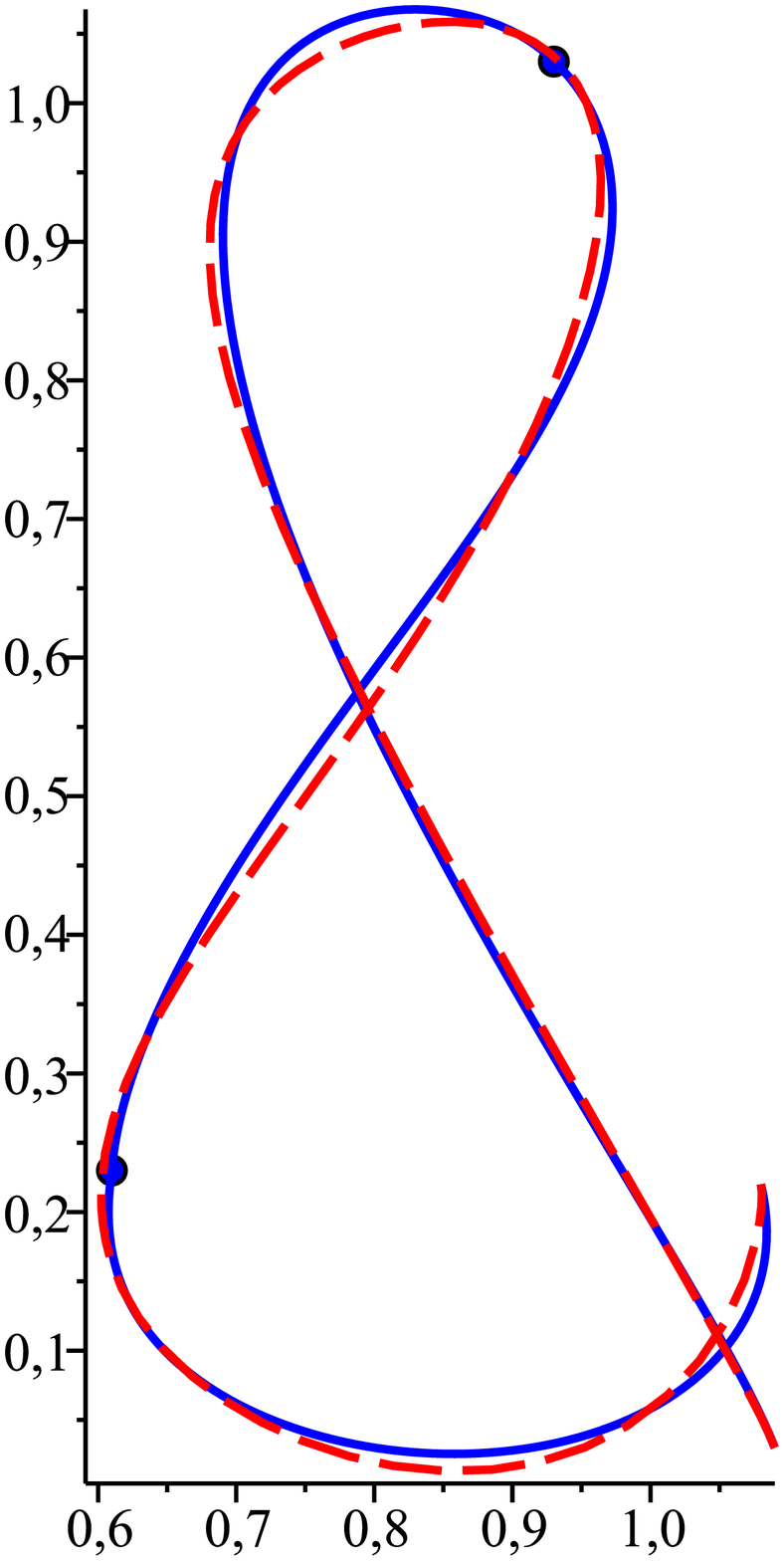}}
\subfloat[]{\label{Fig:1b}\includegraphics[width=0.38\textwidth]{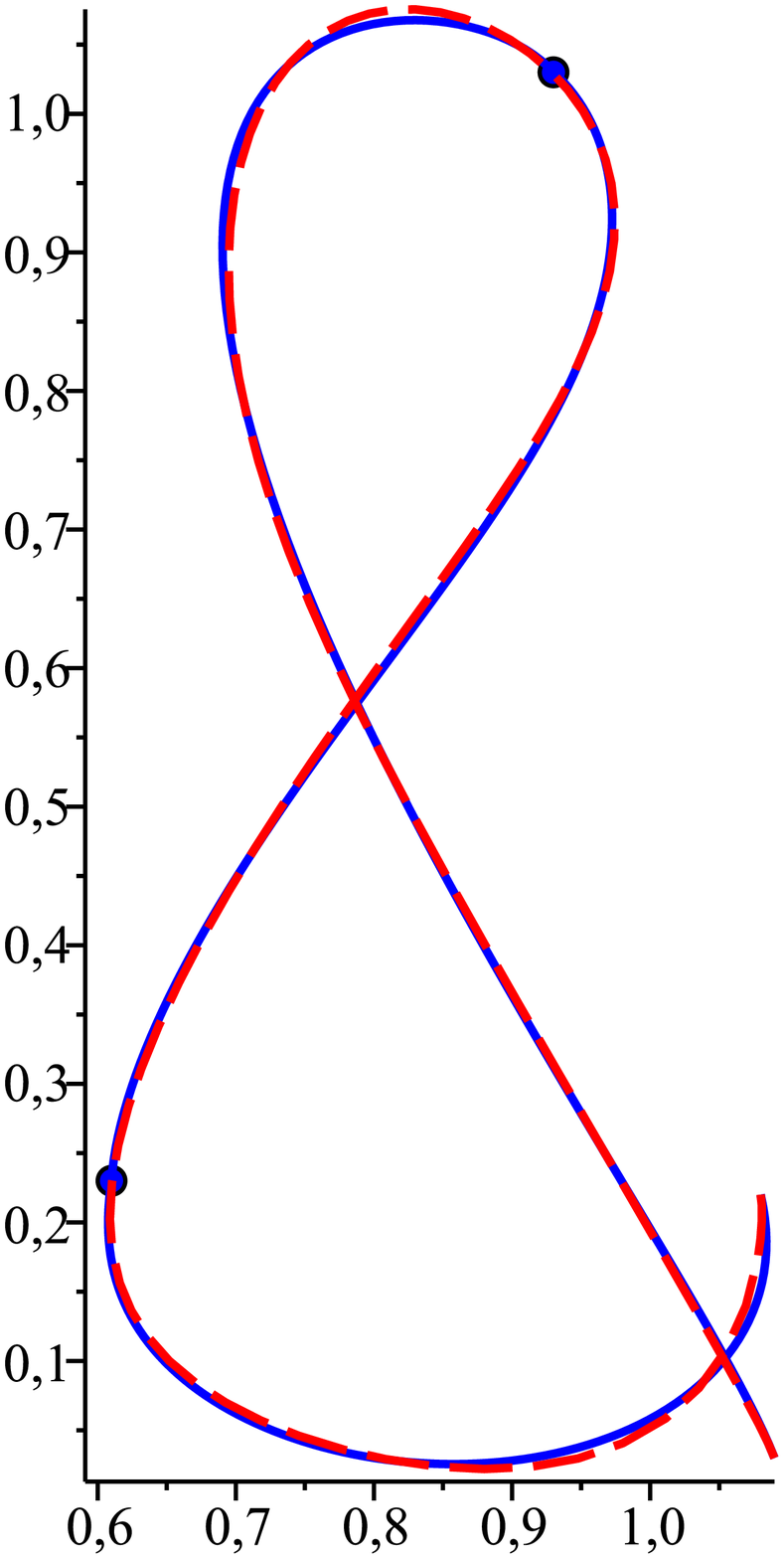}}
\end{tabular}
    \caption{Merging of three  segments of the composite B\'{e}zier curve. Original curve (blue solid line) and merged curve
    (red dashed line) with parameters: (a) $m = 10$, $k = 3$, $l = 2$, and (b) $m = 12$, $k = 3$, $l = 2$.}
\end{center}
\end{figure}

\end{exmp}

\begin{exmp}\label{Ex:Penguin}
The curve ``Penguin'' is formed by two composite B\'ezier curves.
The left curve has four cubic segments, with the control points
$\{(0.31,0.23),$ $(0.35,0.19),$ $(0.39,0.23),$ $(0.37,0.26)\}$,
$\{(0.37,0.26),$ $(0.21,0.54),$ $(0.53,0.77),$ $(0.21,0.76)\}$,
$\{(0.21,0.76),$ $(0.1,0.76),$ $(0.5,$\\$0.88),$ $(0.42,0.79)\}$, and
$\{(0.42,0.79),$ $(0.26,0.76),$ $(0.23,0.92),$ $(0.34,0.94)\}$, respectively.
The  right curve is composed of three cubic segments having control points
$\{(0.34,0.94),$ $(0.74,0.99),$ $(0.67,0.19),$ $(0.56,0.21)\}$,
$\{(0.56,0.21),$ $(0.19,0.32),$ $(0.62,1.05),$ $(0.56,0.61)\}$, and
$\{(0.56,0.61),$ $(0.5,0.24),$ $(0.41,0.41),$ $(0.5,0.64)\}$, respectively.
Formula \eqref{E:part} gives
$t_0 = 0,\ t_1 \doteq 0.08,\ t_2 \doteq 0.55,\ t_3 \doteq 0.78,\ t_4 = 1$ for the left curve,
and  $t_0 = 0,\ t_1 \doteq 0.42,\ t_2 \doteq 0.78,\ t_3 = 1$ for the right one.
Results of separate merging of segments of both  curves  can be seen in Table~\ref{Tab:ex2}.
Two selected cases are shown on Figures~\ref{Fig:2a}~and~\ref{Fig:2b}.
\clearpage

\begin{table}[H]
\captionsetup{margin=0pt, font={scriptsize}}
\ra{1}
\centering
\scalebox{0.88}{
\begin{tabular}{@{}ccccccccccccccccccccc@{}}
\toprule \multicolumn{9}{c}{Left curve} & \phantom{ab} & \multicolumn{9}{c}{Right curve}
\\ \cmidrule{1-9} \cmidrule{11-19}
$m$ & \phantom{a} & $k$ & \phantom{a} & $l$ & \phantom{a} & $E_{2}$ & \phantom{a} & $E_{\infty}$ & \phantom{abc} &
$m$ & \phantom{a} & $k$ & \phantom{a} & $l$ & \phantom{a} & $E_{2}$ & \phantom{a} & $E_{\infty}$
\\ \midrule
$12$ & & $1$ & & $1$ & & $7.45E{-}3$ & &$1.90E{-}2$ & & $10$ & & $1$ & & $1$ & &$1.28E{-}2$ & &$3.51E{-}2$\\
  &   &  $1$ & & $2$ & &$1.05E{-}2$ & &$2.69E{-}2$ & & $$ & & $2$ & & $1$ & &$1.28E{-}2$ & &$3.48E{-}2$\\
  &  &  $2$ &  & $1$ & &$7.85E{-}3$ & &$1.93E{-}2$ & & $$ & & $1$ & & $2$ & &$1.29E{-}2$ & &$3.49E{-}2$\\
  & &   $2$ &  & $2$ & &$1.10E{-}2$ & &$2.85E{-}2$ & & $$ & & $2$ & & $2$ & &$1.30E{-}2$ & &$3.44E{-}2$\\
\midrule
$13$ & & $1$ & & $1$ & &$6.68E{-}3$ & &$1.45E{-}2$ & & $12$ & & $1$ & & $1$ & &$9.01E{-}3$ & &$3.00E{-}2$\\
  &   &  $1$ & & $2$ & &$7.80E{-}3$ & &$1.64E{-}2$ & & $$ & & $2$ & & $1$ & &$1.02E{-}2$ & &$3.27E{-}2$\\
  &  &  $2$ &  & $1$ & &$7.28E{-}3$ & &$1.48E{-}2$ & & $$ & & $1$ & & $2$ & &$1.14E{-}2$ & &$2.98E{-}2$\\
  & &   $2$ &  & $2$ & &$8.53E{-}3$ & &$1.71E{-}2$ & & $$ & & $2$ & & $2$ & &$1.23E{-}2$ & &$3.25E{-}2$\\
  \midrule
$14$ & & $1$ & & $1$ & &$4.39E{-}3$ & &$1.19E{-}2$ & & $13$ & & $1$ & & $1$ & &$8.65E{-}3$ & &$2.83E{-}2$\\
  &   &  $1$ & & $2$ & &$4.51E{-}3$ & &$1.27E{-}2$ & & $$ & & $2$ & & $1$ & &$9.16E{-}3$ & &$2.81E{-}2$\\
  &  &  $2$ &  & $1$ & &$4.86E{-}3$ & &$1.17E{-}2$ & & $$ & & $1$ & & $2$ & &$1.11E{-}2$ & &$2.98E{-}2$\\
  & &   $2$ &  & $2$ & &$5.08E{-}3$ & &$1.30E{-}2$ & &  & & $2$ & & $2$ & &$1.16E{-}2$ & & $2.98E{-}2$
\end{tabular}}
\caption{Least-squares and maximum errors for separate merging of segments of two composite B\'{e}zier curves with constraints.
\label{Tab:ex2}}
\end{table}

\begin{figure}[H]
\captionsetup{margin=0pt, font={scriptsize}}
\begin{center}
\setlength{\tabcolsep}{0mm}
\begin{tabular}{c}
\subfloat[]{\label{Fig:2a}\includegraphics[width=0.51\textwidth]{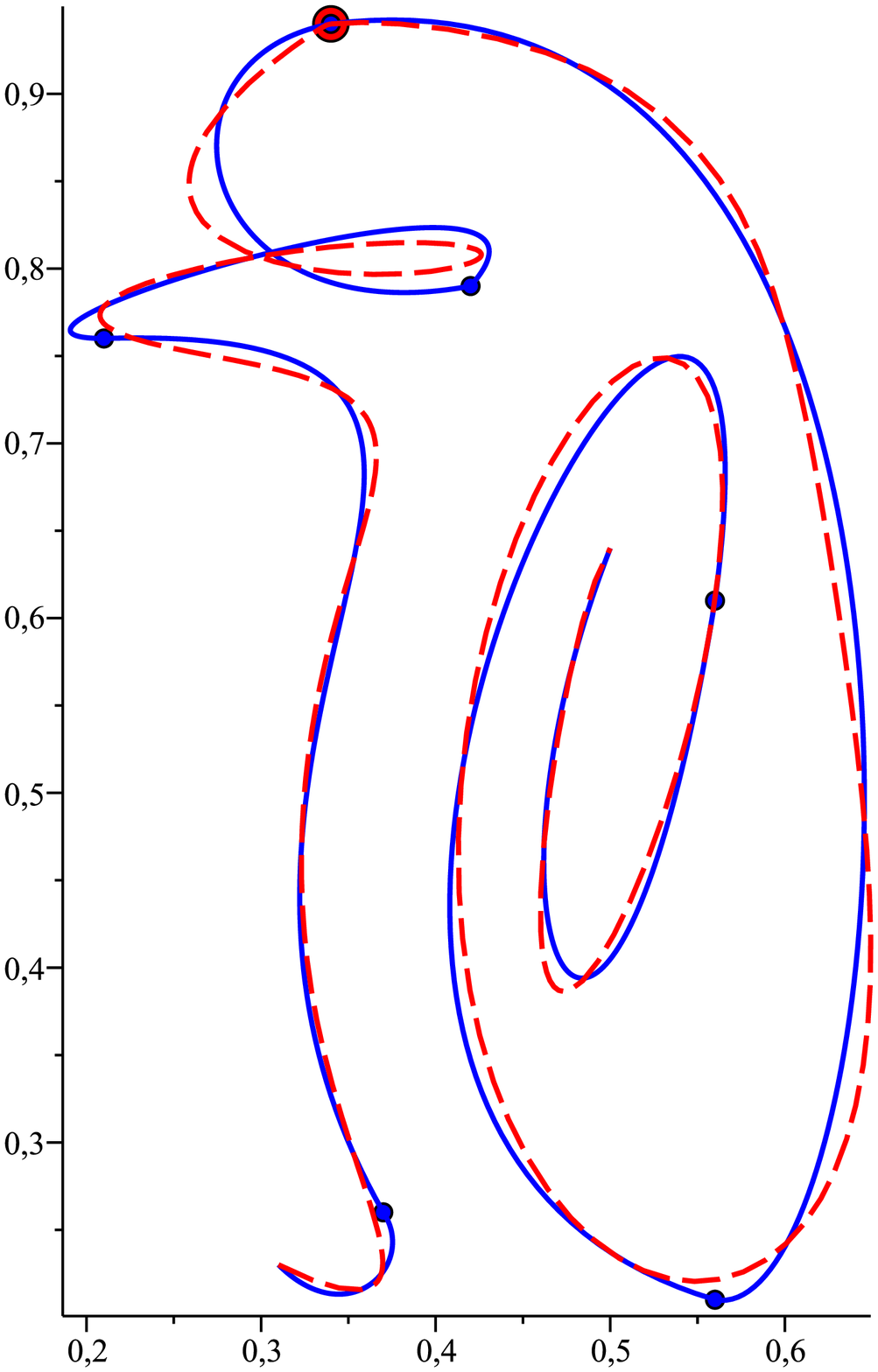}}
\subfloat[]{\label{Fig:2b}\includegraphics[width=0.51\textwidth]{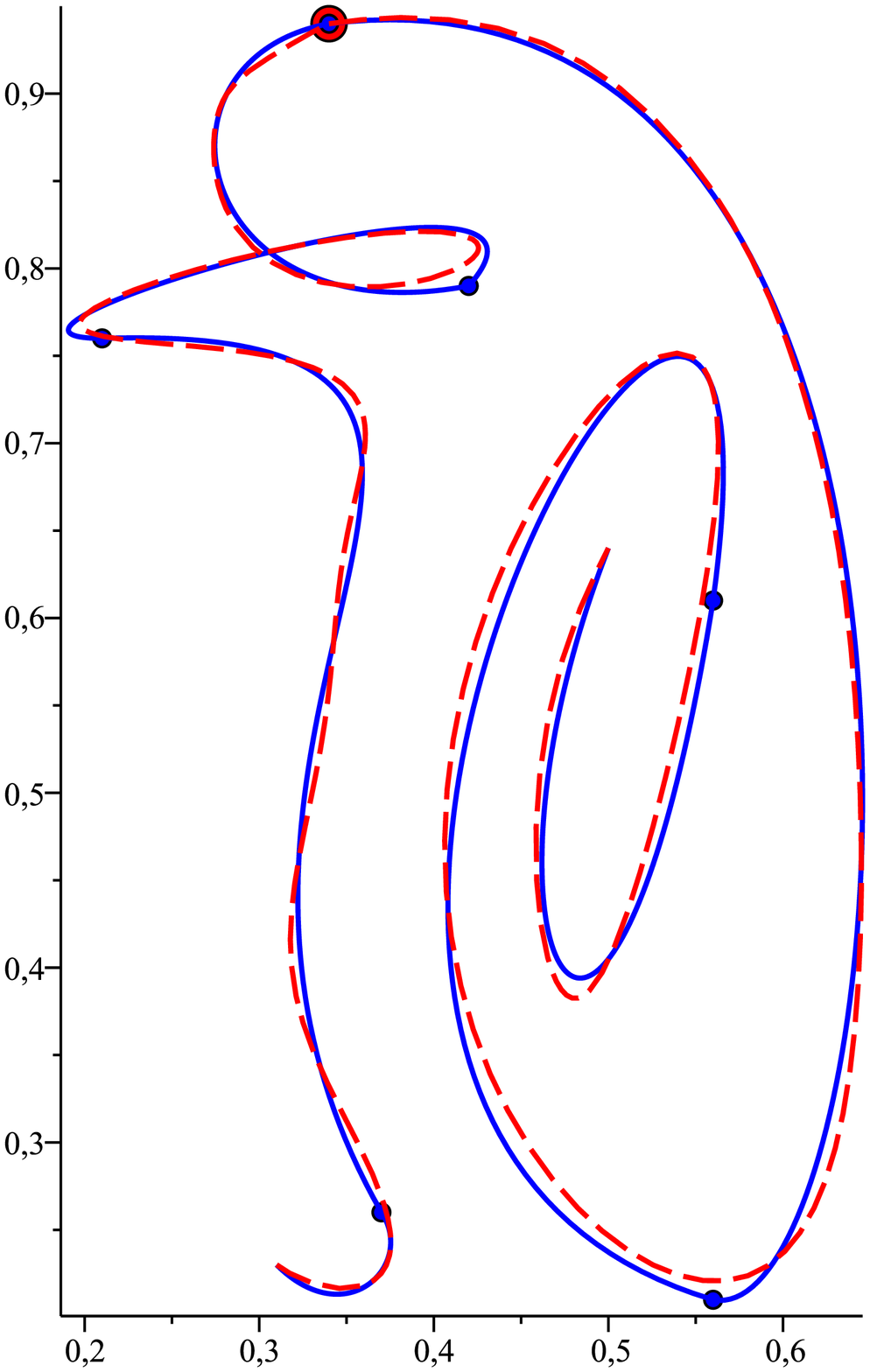}}
\end{tabular}
\caption{Separate merging of segments of two composite B\'{e}zier curves with constraints. Original curves (blue solid line) and
    merged curves (red dashed line).  (a) Left curve: $m = 12$, $k = 1$, $l = 2$;
    right curve: $m = 10$, $k = 2$, $l = 1$. (b) Left curve: $m = 14$, $k = 2$, $l = 2$;
     right curve: $m = 13$, $k = 2$, $l = 2$.}
\end{center}
\end{figure}
\end{exmp}

\clearpage

%%\newpage
%%%%%%%%%%%%%%%%%%%%%%%%%%%%%%%%%%%%%%%%%%%%%%%%%%%%%%%%%%%%%%%%%%%%%%%%%%%%%%%%%%
%%%%%%%%%%%%%%%%%%%%%%%%%%%%%%%%%%%%%%%%%%%%%%%%%%%%%%%%%%%%%%%%%%%%%%%%%%%%%%%%%%
\section{Conclusions}
                    \label{S:Concl}%%SSSSSSSSSSSSSSSSSSSSSSSSSSSSSSSSSSSSSSSSSSS
%%%%%%%%%%%%%%%%%%%%%%%%%%%%%%%%%%%%%%%%%%%%%%%%%%%%%%%%%%%%%%%%%%%%%%%%%%%%%%%%%%
%%%%%%%%%%%%%%%%%%%%%%%%%%%%%%%%%%%%%%%%%%%%%%%%%%%%%%%%%%%%%%%%%%%%%%%%%%%%%%%%%%

We have proposed a novel approach to the problem
of merging  of   multiple adjacent B\'ezier curves, with the endpoints continuity constraints.
\begin{comment}
We have shown that it is possible to generalize dual Bernstein
polynomials approach to compute the control points of the merged curve.

Thanks
to  using fast schemes of evaluation of  certain connections
involving  Bernstein and dual
Bernstein polynomials,  the complexity of our algorithm is $O(sm^2)$,
which should be compared to the complexity $O(sm^3)$ of the only existing multiple merging method
\cite{CW08}.
\end{comment}
We have shown that, contrary to some earlier opinions \cite{Lu15}, it is possible to generalize dual Bernstein
polynomials approach to compute the control points of the merged curve.
Thanks
to  using fast schemes of evaluation of  certain connections
involving  Bernstein and dual
Bernstein polynomials,  the complexity of our algorithm is $O(sm^2)$,
which should be compared to the complexity $O(sm^3)$ of the existing multiple merging methods
\cite{CW08,Lu15}.

As for our future work, we plan to study the above merging problem with $G^{k,l}$ continuity constraints.

%%%%%%%%%%%%%%%%%%%%%%%%%%%%%%%%%%%%%%%%%%%%%%%%%%%%%%%%%%%%%%%%%%%%%%%%%%%%%
%%%%%%%%%%%%%%%%%%%%%%%%%%%%%%%%%%%%%%%%%%%%%%%%%%%%%%%%%%%%%%%%%%%%%%%%%%%
%%%%%%%%%%%%%%%%%%%%%%%%%%%%%%%%%%%%%%%%%%%%%%%%%%%%%%%%%%%%%%%%%%%%%%%%%%%


\begin{thebibliography}{99}
\itemsep2pt

\bibitem{CW08} M. Cheng, G. Wang, Approximate merging of multiple B\'ezier segments,
	       Progress in Natural Science 18 (2008), 757--762.

\bibitem{Far02} G.\,E. Farin,
                Curves and Surfaces for Computer-Aided
                Geometric Design. A Practical Guide, fifth edition,
                Academic Press,  Boston, 2002.

\bibitem{Hos87} J. Hoschek, Approximate conversion of spline curves,
		Computer Aided Geometric Design 4 (1987), 59--66.

\bibitem{HTJS01} S. Hu, R. Tong, T. Ju, J. Sun, Approximate merging of a pair of B\'ezier curves,
 		 Computer-Aided Design 33 (2001), 125--136.

 		 \bibitem{Jue98} B. J\"uttler, The dual basis functions for the Bernstein polynomials,
		Advances in Computational Mathematics 8 (1998),  345--352.
		
\bibitem{LW11a}  S. Lewanowicz,  P. Wo\'zny,   Multi-degree reduction of tensor product
		B\'ezier surfaces with general constraints,
		Applied Mathematics and Computation  217 (2011), 4596--4611.

\bibitem{LW11b}   S. Lewanowicz,  P. Wo\'zny,
                B\'ezier representation of the constrained dual Bernstein polynomials,
                Applied Mathematics and Computation 218 (2011), 4580--4586.       		

\bibitem{Lu14}  L. Lu, An explicit method for $G^3$ merging of two B\'ezier curves,
		Journal of Computational and  Applied  Mathematics 260 (2014), 421--433.
		
\bibitem{Lu15} L. Lu, Explicit algorithms for multiwise merging of  B\'ezier curves,
                Journal of Computational and  Applied  Mathematics 278 (2015), 138--148.		

\bibitem{Lu13} L. Lu, Effective $C^1G^2$-merging of Two B\'ezier Curves by Matrix Computation,
               International Journal of Advancements in Computing Technology 5 (2013), 1117--1123.

\bibitem{THH03} C. Tai, S. Hu, Q. Huang, Approximate merging of B-spline curves via knot adjustment and
		constrained optimization, Computer-Aided Design 35 (2003), 893--899.
		
\bibitem{WL09}  P. Wo\'zny,  S. Lewanowicz,  Multi-degree reduction of B\'ezier curves
		with constraints, using dual Bernstein basis polynomials,
		Computer Aided Geometric Design 26 (2009), 566--579.
		
\bibitem{ZW09}  P. Zhu, G. Wang, Optimal approximate merging of a pair
		of B\'ezier curves with $G^2$-continuity,
		Journal of Zhejiang University SCIENCE A 10 (2009), 554--561.		
\end{thebibliography}
\end{document}